\documentclass [10pt,a4paper]{article}

\usepackage{amssymb}
\usepackage{enumerate}
\usepackage{amsmath}
\usepackage{amsfonts,mathrsfs}
\usepackage{amsthm}
\usepackage{graphicx}

\usepackage[applemac]{inputenc}
\usepackage[width=1.1\textwidth]{caption}
\usepackage[usenames,dvipsnames]{color}
\usepackage[T1]{fontenc}
\usepackage{lpic}

\usepackage{color}
\definecolor{green}{rgb}{0,0.8,0.5}

\usepackage{hyperref}
\hypersetup{backref,colorlinks=true,linkcolor=blue,citecolor=red}
\newcounter{teoremaganso}

\newenvironment{trivlistparm}[2]{\trivlist
\item[{#1 #2}]}{\endtrivlist}

\newenvironment{prooftext}[1]{\trivlistparm{\bfseries}{#1}}{\Qed\endtrivlistparm}
\newenvironment{prova}{\trivlistparm{\bfseries}{Proof.}}{\Qed\endtrivlistparm}

\catcode`\@=11

\def\resetthefootnote{\renewcommand{\thefootnote}{\@arabic\c@footnote} }
\def\@principiremex#1{\trivlist
 \item[\hskip \labelsep{\bfseries #1\ \thetheo.}]\ignorespaces}
\def\opar@principiremex#1[#2]{\trivlist
 \item[\hskip \labelsep{\bfseries #1\ \thetheo\ (#2).}]\ignorespaces}

\newcommand{\newTHEOremrom}[2]{\newenvironment{#1}{\refstepcounter{theo}\@ifnextchar[{\opar@principiremex{#2}}
{\@principiremex{#2}}}{\qedB\endtrivlist}} \catcode`\@=12
\DeclareMathSymbol{\square}{\mathord}{AMSa}{"03}
\newcommand{\qedB}{\nopagebreak\hspace*{\fill}$\square$\par}
\newcommand{\Qed}{\nopagebreak\hspace*{\fill}{\vrule width6pt height6pt depth0pt}\par}

\newtheorem {theo} {Theorem} [section]
\newtheorem {prop} [theo] {Proposition}

\newtheorem {lem} [theo] {Lemma}
\newtheorem {bigtheo} [teoremaganso] {Theorem}

\newTHEOremrom {defi} {Definition}
\newTHEOremrom {obs} {Remark}
\newTHEOremrom {ex} {Example}

\newcommand{\refc}[1]{\mbox{$(\ref{#1})$}}
\newcommand{\secc}[1]{Section~\ref{#1}}
\newcommand{\teoc}[1]{Theorem~\ref{#1}}
\newcommand{\propc}[1]{Proposition~\ref{#1}}

\newcommand{\lemc}[1]{Lemma~\ref{#1}}
\newcommand{\defic}[1]{Definition~\ref{#1}}
\newcommand{\obsc}[1]{Remark~\ref{#1}}

\newcommand{\figc}[1]{Figure~\ref{#1}}

\newcommand{\C}{\mathbb{C}}

\newcommand{\R}{\mathbb{R}}
\newcommand{\Z}{\mathbb{Z}}
\newcommand{\N}{\mathbb{N}}

\newcommand{\vp}{\varphi}

\newcommand{\ka}{\kappa}

\newcommand{\be}{\vartheta}
\newcommand{\g}{\gamma}
\newcommand{\la}{\lambda}
\newcommand{\A}{a}
\newcommand{\Am}{a_{M}}
\newcommand{\rr}{r}

\newcommand{\De}{\Delta}

\newcommand{\infa}{K}

\renewcommand{\P}{\mathscr{P}}

\newcommand{\Bs}{\mathscr{B}_{\sigma}}
\newcommand{\bb}{\bar{\be}}

\newcommand{\ii}{\ensuremath{(x_{\ell},x_r)}}

\newcommand{\op}{\ensuremath{\mbox{\rm o}}}

\def\map#1#2#3{\mbox{${#1}\!:{#2}\longrightarrow{#3}$}}

\newcommand{\sist}[2]{
  \left\{\!
   \begin{array}{l}
    \dot x=#1 \\[2pt] \dot y=#2
   \end{array}
  \right.
}

\title{\textbf{On the wave length of smooth periodic \\ traveling waves  of the Camassa-Holm equation}
\footnotetext{A.~Geyer is  supported by the FWF project J3452 ``Dynamical
Systems Methods in Hydrodynamics'' of the Austrian Science Fund. J.~Villadelprat is partially supported by the MEC/FEDER grant MTM2008-03437.}}

\author{A. Geyer and J. Villadelprat
\\*[.1truecm]
{\small \textsl{Departament de Matem\`atiques}}
\\*[-.05truecm]
{\small \textsl{Universitat Aut\`onoma de Barcelona, Barcelona,
Spain}}
\\*[.1truecm]
{\small \textsl{Departament d'Enginyeria Inform\`{a}tica i Matem\`{a}tiques}}
\\*[-.05truecm]
{\small \textsl{Universitat Rovira i Virgili, Tarragona, Spain}}}

\date{}

\begin{document}
\maketitle
\begin{abstract}
This paper is concerned with the wave length $\lambda$ of smooth periodic traveling wave solutions of the Camassa-Holm equation. The set of these solutions can be parametrized using the wave height $a$ (or ``peak-to-peak amplitude''). Our main result establishes  monotonicity properties of the map $a\longmapsto \lambda(a)$, i.e.,
the wave length as a function of the wave height. We obtain the explicit bifurcation values, in terms of the parameters associated to the equation, which distinguish between the two possible qualitative behaviours of $\lambda(a)$, namely monotonicity and unimodality. 
The key point is to relate $\lambda(a)$ to the period function of a planar differential system with a quadratic-like first integral, and to apply a criterion which bounds the number of critical periods for this type of systems.
\end{abstract}

%%%%%%%%%%%%%%%%%%%%%%%%%%%%%%%%%%%%%%%%%%%%%%%%%%%%%%%%%%%%%%%%%%%%%%%%%%

\section{Introduction and main result}
\label{sect: Intro}
The Camassa-Holm (CH) equation 
\begin{equation}
 \label{CH}
  u_t + 2\ka\,u_x -u_{txx} + 3\, u\,u_x  = 2\,u_x u_{xx} + u\,u_{xxx},
       \qquad x \in \R, \; t>0,
 \end{equation}
 arises as a shallow water approximation of the Euler equations for inviscid, incompressible and homogenous fluids propagating over a flat bottom, where $u(x,t)$ describes the horizontal velocity component and $\ka\in\R$ is a parameter related to the critical shallow water speed. 
This equation was first derived by Fokas and Fuchssteiner \cite{Fuchssteiner1981} as an abstract bi-Hamiltonian equation with infinitely many conservation laws, and later re-derived by Camassa and Holm  \cite{Camassa1993} from physical principles. 
For a discussion on the relevance and applicability of the CH equation in the context of water waves we refer the reader to Johnson \cite{Joh02,Joh03b, Joh03a} and more recently Constantin and Lannes \cite{ConLan09}. 
We point out that for a large class of initial conditions the CH equation is an integrable infinite-dimensional Hamiltonian system \cite{Teschl2009,Con01,Constantin2006c,  ConLen2003,Constantin2002,Joh03b},
and it is known that the solitary waves of CH are solitons which are orbitally stable \cite{
Constantin2002, DikaMoli2007}. 
Some classical solutions of the CH equation develop singularities in finite time in the form of wave breaking: 
the solution remains bounded but its slope becomes unbounded \cite{Camassa94,Constantin1998g, Constantin1998a,Constantin1998,Dan01,Li2000,McKean1998}.
After blow-up the solutions can be recovered in the sense of global weak solutions, see~\cite{Bressan2007a,Bressan2007} and also \cite{Holden2007, Raynaud2012}.

\begin{figure}[t]
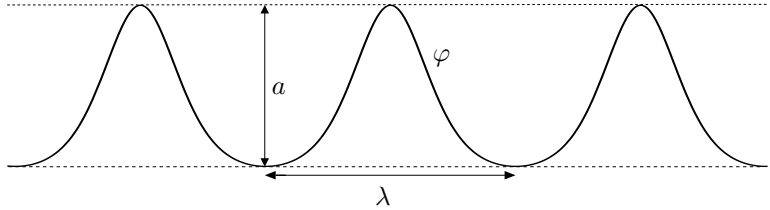

  \centering    
 \begin{lpic}[l(2mm),r(2mm),t(2mm),b(2mm)]{Figure0(0.4)}
  \lbl[l]{140,50;$\vp$}
  \lbl[l]{87,40;$a$}
  \lbl[l]{121,4.5;$\la$}
  
 \end{lpic}     
\caption{Smooth periodic TWS $\vp$ of CH with wave length $\la$ and wave height $a$.}
\label{Fig-wave}
\end{figure}

In the present paper, we consider traveling wave solutions of the form 
\begin{equation}
 \label{e-Ansatz}
    u(x,t)=\vp(x-c\,t),
\end{equation}
for $c\in\R$ and some function $\vp: \R\rightarrow \R$. We denote $s = x- ct$ the independent variable in the moving frame. Inserting the  Ansatz \eqref{e-Ansatz} into equation \eqref{CH} and integrating once we obtain the corresponding equation for traveling waves,
\begin{equation}
\label{e-CHODE}
 \vp''(\vp-c) + \frac{(\vp')^2}{2} + \rr  + (c-2\ka)\,\vp - \frac{3}{2}\vp^2 = 0,
\end{equation} 
where $\rr  \in \R$ is a constant of integration and the prime denotes derivation with respect to~$s$. 
A solution~$\vp$ of~\eqref{e-CHODE} is called a \emph{traveling wave solution} (TWS) of the Camassa-Holm equation~\eqref{CH}. Lenells~\cite{Lenells2005a} provides a complete classification of all (weak) traveling wave solutions of the Camassa-Holm equation. In the present paper, we focus on \emph{smooth periodic} TWS of the Camassa-Holm equation, which can be shown to have a unique maximum and minimum per period, see~\cite{Lenells2005a}.
In the context of fluid dynamics the period of such a solution is called \emph{wave length}, which we will denote by $\lambda$. The difference between the maximum (wave crest) and the minimum (wave trough) is called \emph{wave height}, see~Figure~\ref{Fig-wave}, which we will denote by $a$ (in some contexts this quantity is also called ``peak-to-peak amplitude'').

The aim of this paper is to study the dependence of the wave length $\lambda$ of smooth periodic TWS of the Camassa-Holm equation \eqref{CH} on their wave height $a$. Our main result shows that $\lambda(a)$ is a well-defined function and that it is either monotonous or unimodal. More precisely:

\begin{bigtheo}\label{thm}
Given $c, \ka$ with $c\neq -\ka$, there exist real numbers $\rr_1<\rr_{b_1} < \rr_{b_2} < \rr_2$ such that the differential equation \eqref{CH} has smooth periodic TWS of the form \eqref{e-Ansatz} if, and only if, the integration constant $r$ in \eqref{e-CHODE} belongs to the interval $(\rr_1,\rr_2)$. For such $\rr\in (\rr_1,\rr_2)$, the set of  smooth periodic TWS form a continous family $\{\vp_\A\}_{\A\in(0,\Am)}$ parametrized by the wave height~$\A$. Furthermore, the wave length $\la= \la(\A)$ of $\vp_{\A}$ satisfies the following:
\begin{enumerate}[$(a)$]
 \item  If $\rr\in(\rr_1,\rr_{b_1}]$, then $\la(a)$ is monotonous increasing.
 \item If $\rr\in(\rr_{b_1},\rr_{b_2})$, then $\la(a)$ has a unique critical point which is a maximum.
  \item If $\rr\in[\rr_{b_2},\rr_2)$, then $\la(a)$ is monotonous decreasing.
\end{enumerate}
Finally, these are the only possible scenarios for smooth periodic TWS of the CH  equation.
\end{bigtheo}

We point out, see \propc{p-prop1}, that if $c=-\ka$ then there are no smooth periodic TWS of the form~\eqref{e-Ansatz}. The exact values of the bifurcation parameters  $\rr_1, \rr_{b_1}, \rr_{b_2}$ and $\rr_2$ in terms of $c$ and $\ka$ can be found in the proof of \teoc{thm} at the end of \secc{s-pre}. In this regard, we remark that the expressions $\rr_1, \rr_{b_1}$ and~$\rr_2$ also appear in~\cite[p.~402]{Lenells2005a}, but they serve as bifurcation values for a different type of property: they  define the boundaries of parameter regions  where the various types of weak TWS (smooth, peaked or cusped waves, \ldots) can occur. 
It should also be observed that a description on how the wave length of TWS of CH depends on parameters  may be found in the last section of \cite{Lenells2005a}, where level sets of TWS with the same wave length are described. Furthermore,  it is shown that there exist peakons and cuspons with arbitrarily small wave length.  In contrast, we will show that for smooth periodic TWS the wave length cannot be arbitrarily small, see \obsc{remark}.

The paper is organized as follows. In \secc{s-pre} we establish a correspondence between smooth periodic TWS of \eqref{CH} and periodic orbits around the center of a planar differential  system with a quadratic-like first integral, see~\propc{p-prop1}. We observe that the wave length of a smooth periodic solution of \eqref{e-CHODE} is equal to the period of the corresponding periodic orbit. Moreover,  there exists an analytic diffeomorphism which relates the  wave height of a solution of \eqref{e-CHODE} to the energy level of the first integral at the corresponding periodic orbit of the planar system, see~\lemc{leng}. In \teoc{res_periode}, we state the monotonicity properties of the period function of the center of this planar system, which imply \teoc{thm}. 
The proof of \teoc{res_periode} is carried out in \secc{s-mon}. It relies on a result proved in \cite{Vil2014}, which provides a criterion to bound the number of critical periods for this kind of systems.

%%%%%%%%%%%%%%%

\section{Smooth periodic TWS of the Camassa-Holm equation}\label{s-pre}

TWS of the form \eqref{e-Ansatz} of the Camassa-Holm equation \eqref{CH}  correspond to solutions of the equation~\eqref{e-CHODE}.  The next result establishes a correspondence between the smooth periodic solutions of \eqref{e-CHODE} and periodic orbits around the center of an associated planar system. Moreover, it provides a necessary and sufficient condition for the existence of such a center. To this end, recall that the largest punctured neighbourhood of  a  center which consists entirely of periodic orbits  is called \emph{period annulus}, see~\cite{Chicone}.

\begin{prop}\label{p-prop1}
The following holds:
\begin{enumerate}[$(a)$]
\item $\vp$ is a smooth periodic solution of equation \refc{e-CHODE} if, and only if, $(w,v)=(\vp-c,\vp')$ 
         is a periodic orbit of the planar differential system 
\begin{equation}\label{e-sys_w}
       \left\{
      \begin{array}{l}
    w' =v,\\[2pt]
    v' =-\dfrac{F'(w)+\frac{1}{2}\,v^2}{w},
      \end{array}\right.
  \end{equation}
         where 
\begin{equation}\label{e-F}
    \mbox{$F(w)\!:= \alpha w + \beta  w^2 -\frac{1}{2} w^3,$ with $\alpha \!:=r-2\ka c- \frac{1}{2} c^2$ and $\beta \!:= - (c+\ka)$.}
\end{equation}
\item The function $\hat H(w,v)\!:=\frac{1}{2}w v^2+F(w)$ is a first integral of the differential system~\refc{e-sys_w}.
\item Every periodic orbit of system \refc{e-sys_w} belongs to the period annulus $\mathscr P$ of a center, 
         which exists if, and only if, $-2\beta ^2 < 3\alpha  < 0$ is verified. 
 \end{enumerate}
\end{prop}

\begin{prova} 
The assertion in $(b)$ is straightforward. In order to prove $(a)$ we first note that \refc{e-CHODE} can be written as $\vp''(\vp-c)+\frac{1}{2}(\vp')^2+F'(\vp-c)=0,$ where $F$ is defined in \eqref{e-F}. Accordingly, $\vp$ is a solution of~\refc{e-CHODE} with $\vp(s)\neq c$ for all $s$ if, and only if, $s\longmapsto (w,v)=\bigl(\vp-c,\vp'\bigr)(s)$ is a solution of the differential system~\refc{e-sys_w}. We claim that $\vp(s)\neq c$ for all~$s\in\R$ in case that~$\vp$ is smooth and periodic, i.e.~$\vp(s+T)=\vp(s)$ for some $T>0.$ Clearly, $(a)$ will follow once we show the claim. With this aim in view note  that if $\vp$ is a smooth periodic solution of \eqref{e-CHODE} then the set $\mathscr C\!:=\bigl\{(w,v)=\bigl(\vp-c,\vp'\bigr)(s);s\in\R\bigr\}$ describes a smooth loop. We will show that $\mathscr C$ cannot intersect $\{w=0\}$. We can rule out that $\vp\equiv c$ because a constant function is not periodic. Hence suppose that there exist $s_0$ and $s_1$ such that $\vp(s)\neq c$ for all $s\in (s_0,s_1)$ and $\vp(s_1)=c$. Then, for $s\in (s_0,s_1),$ $(w,v)=\bigl(\vp-c,\vp'\bigr)(s)$ is a solution of the differential system \refc{e-sys_w} that tends to the point $p_1\!:=\bigl(0,\vp'(s_1)\bigr)$ as $s\longrightarrow s_1.$ Since $\hat H(p_1)=0$ and by the continuity of $\hat H,$ it turns out that $\mathscr C$ is inside the zero level set of $\hat H$. An easy computation shows that $\hat H(w,v)=0$ if, and only if, $w=0$ or $(w-\beta )^2-v^2=\beta ^2+2\alpha .$ The second equality describes a hyperbola which intersects $\{w=0\}$ if, and only if, $\alpha \leqslant 0.$ In any case it is not possible that $\{\hat H=0\}$ contains a smooth loop. So the claim is true and $(a)$ follows. 
\begin{figure}[t]
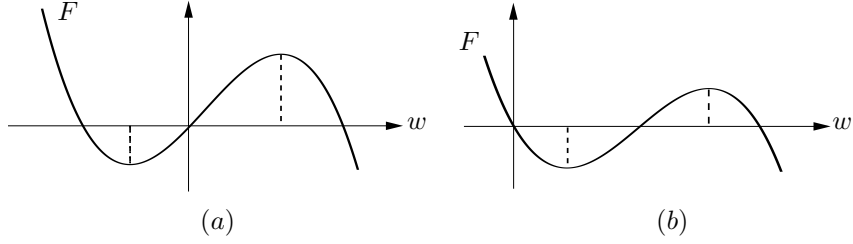

  \centering    
 \begin{lpic}[l(2mm),r(2mm),t(0mm),b(5mm)]{Figure1(0.4)}
      \lbl[l]{80,-10;$(a)$}
      \lbl[l]{230,-10;$(b)$}
      \lbl[l]{33,60;$F$}
      \lbl[l]{165,50;$F$}
      \lbl[l]{148,23;$w$}
      \lbl[l]{288,23;$w$}     
 \end{lpic}     
\caption{Sketch of the graph of $F$: $(a)$ when $\alpha >0$, $(b)$ when $\alpha <0$ and $\beta >0$.}
\label{Fig-center}
\end{figure}

In order to show $(c)$ recall that the differential system~\refc{e-sys_w} has a  first integral, and consequently there are no limit cycles and the periodic orbits form period annuli. A periodic orbit must surround at least one critical point of the differential system, which are of the form $(w,v)=(\hat w,0)$ with $\hat w\neq 0$ and $F'(\hat w)=0.$ The determinant of the Jacobian of the vector field at such a point is $\det J_{(\hat w,0)}=\frac{F''(\hat w)}{\hat w}.$ A straightforward computation shows that $F'(w)=0$ if, and only if, $w=\frac{2\beta  \pm \sqrt{4\beta ^2+6\alpha }}{3}.$ If $\alpha >0$, then~$F$ has a minimum on $w<0$ and a maximum on $w>0$ (see \figc{Fig-center}), 
which both correspond to saddle points of system~\eqref{e-sys_w}. 
Thus, by applying the Poincar\'{e}-Bendixon Theorem (see for instance~\cite{Sotomayor1979}), no periodic orbit is possible in case that $\alpha >0$.
Similarly, if $\alpha =0$ and $\beta \neq 0,$ then there is only one critical point, which is a saddle, whereas if $\alpha =0$ and $\beta = 0,$ then there are no critical points. Hence there are no periodic orbits in case that $\alpha =0.$ 
Finally let us discuss the case $\alpha <0.$ If $\alpha <0$ and $2\beta ^2+3\alpha <0$, then there are no critical points, which prevents the differential system from having periodic orbits. If $\alpha <0$ and $2\beta ^2+3\alpha =0$, then there exits a unique critical point which is a cusp, and can not be surrounded by a periodic orbit. 
If $\alpha <0$ and $2\beta ^2+3\alpha >0$ then, see \figc{Fig-center}, $F$ has its two local extrema, which are located on $w<0$ in case that $\beta <0$, and on $w>0$ in case that $\beta >0.$ In both cases one extremum yields  a saddle and the other a center of system~\eqref{e-sys_w}. By applying the Poincar\'{e}-Bendixon Theorem one can easily conclude that the set of periodic orbits forms a punctured neighbourhood of the center, and that no other period annulus is possible. This proves $(c)$. 
\end{prova}

It is now necessary to introduce some notation. 
\begin{defi}
\label{Def}
Let $\vp$ be a smooth periodic solution of the differential equation~\refc{e-CHODE}.  We denote by $a_{\vp}$ its wave height. By \propc{p-prop1}, 
$(w,v)=(\vp-c,\vp')$ is a periodic orbit inside the period annulus $\P$ of the differential system  \eqref{e-sys_w}, which we denote by $\gamma_{\vp}$. Since $\hat H$ is a first integral of system \eqref{e-sys_w}, the orbit $\gamma_{\vp}$ is inside some level curve of $\hat H$, and we denote its energy level by $h_{\vp}$. In addition, let  the center of \eqref{e-sys_w}  be inside the level curve $\{\hat{H}=h_0\}$ and suppose that $\hat{H}(\P)=(h_0,h_1)$. Then $h_{\vp}\in (h_0,h_1)$.
\end{defi}

The following result establishes a relation between the wave height of a smooth periodic solution of  \eqref{e-CHODE} and the energy level of the corresponding periodic orbit of \eqref{e-sys_w}.

\begin{lem}\label{leng}
Suppose that the set $\{\vp\}$ of smooth periodic solutions  of \eqref{e-CHODE} is nonempty. With the notation introduced in \defic{Def}, the following holds:
\begin{enumerate}[$(a)$]
\item The period of the periodic orbit $\gamma_{\vp}$ is equal to the wave length of $\vp.$
\item There exists an analytic diffeomorphism \map{\ell}{(h_0,h_1)}{(0,\Am)} verifying that $\ell(h_{\vp})=a_{\vp}$ 
         for all $\vp$. In addition,~$\ell$ can be analytically extended to $h=h_0$ by setting $\ell(h_0)=0.$

\end{enumerate}
\end{lem}

\begin{prova}
The assertion in $(a)$ is clear. The key point to prove $(b)$ is that the length of the projection of the periodic orbit $\gamma_{\vp}$ on the $w$-axis is $a_{\vp}.$ In order to compute it let us fix that the center of the differential system \refc{e-sys_w} is at the point $(w_c,0)$. Let $(w_{\ell},w_r)$ be the projection of its period annulus $\mathscr P$ on the $w$-axis. Thus $w_{\ell}<w_c<w_r$ and $F'(w)\neq 0$ for all $w\in (w_{\ell},w_r)\setminus\{w_c\},$ whereas $F'(w_c)=0$ and $F''(w_c)\neq 0.$ Then there exits an analytic diffeomorphism $G$ on $(w_{\ell},w_r)$ such that $F(w)=h_0+G(w)^2$, where $h_0=F(w_c).$ Recall that, by definition, the periodic orbit $\gamma_{\vp}$ is inside the energy level $\{\hat H=h_{\vp}\}$. Since $\hat H(w,0)=h$ if and only if $h_0+G(w)^2=h,$ we have that $\gamma_{\vp}$ intersects the $w$-axis at the points $p_{\pm}(\vp)=\bigl(G^{-1}(\pm\sqrt{h_{\vp}-h_0}\,),0\bigl)$. Hence the length of its projection on the $w$-axis is $a_{\vp}=\ell(h_{\vp})\!:=G^{-1}(\sqrt{h_{\vp}-h_0}\,)-G^{-1}(-\sqrt{h_{\vp}-h_0}\,)$. A  straightforward argument shows that $\ell$ is an analytic diffeomorphism on $(h_0,h_1)$ and that it can be analytically extended to $h=h_0$ setting $\ell(h_0)=0.$ This shows $(b)$ and completes the proof. 
\end{prova}

\begin{obs}\label{equivalencia}
It is clear that the energy levels of $\hat H$  parameterize the set of periodic orbits inside $\P$. Thus, the set of periodic orbits of \refc{e-sys_w} forms a  continuous family $\{\gamma_{h}\}_{h\in (h_0,h_1)}$. Consequently, and thanks to \propc{p-prop1} and \lemc{leng}, we can assert that the set of smooth periodic solutions of~\refc{e-CHODE} forms a continuous family $\{\vp_a\}_{a\in (0,\Am)}$ parameterized by their wave height. We can thus consider the function \map{\lambda}{(0,\Am)}{\R^+} which assigns to each $a\in (0,\Am)$ the wave length of the unique smooth periodic solution of~\refc{e-CHODE} with wave height $a$. \teoc{thm} is concerned precisely with the qualitative properties of this function. We stress that \emph{a priori} it is defined on the set of smooth periodic solutions of~\refc{e-CHODE} rather than on the interval $(0,\Am).$ On account of \lemc{leng}, the wave length $\lambda(a)$ is equal to the period of the periodic orbit of~\refc{e-sys_w} inside the level curve $\{\hat H=\ell^{-1}(a)\}.$ This is the key point in proving \teoc{thm}, as it allows us to deduce qualitative properties of the function $\lambda$ from those of the period function of the center of~\eqref{e-sys_w}. 
\end{obs}

The following technical result, which will be proved in \secc{s-mon},  provides a detailed account on the monotonicity properties of the period function of the center at the origin of the differential system~\eqref{e-sys_w}.
 
\begin{theo}\label{res_periode}
Consider system~\refc{e-sys_w} with $-2\beta ^2 < 3\alpha  < 0$ and define $\be\!:=\frac{1}{6}\!\left(\frac{2|\beta |}{\sqrt{4\beta ^2+6\alpha }}-1\right).$ Then $\be>0$ and the period function of the center of system \eqref{e-sys_w} verifies the following:
 \begin{enumerate}[$(a)$] 
  \item It is monotonous decreasing in case that $\be \in\left(0,-\frac{1}{10}+\frac{1}{15}\sqrt{6}\,\right]$.
  \item It has a unique critical period, which is a maximum, in case that 
           $\be \in\left(-\frac{1}{10}+\frac{1}{15}\sqrt{6},\frac{1}{6}\right)$.
  \item It is monotonous increasing in case that $\be \geqslant\frac{1}{6}$.         
 \end{enumerate}
\end{theo}

We are now in position to prove the main result of the paper.

 \begin{prooftext}{Proof of \teoc{thm}.}
Consider the differential equation \refc{e-CHODE} and define $\alpha =r-2\ka c- \frac{1}{2} c^2$ and $\beta =-(c+\ka)$. It follows from \propc{p-prop1} that the Camassa-Holm equation~\refc{CH} has smooth periodic TWS if, and only if, $-2\beta ^2<3\alpha <0.$ It is easy to see that in terms of the ``intrinsic'' parameters $\ka$ and $c$, these conditions are equivalent to requiring that the integration constant $r$ belongs to the interval $(r_1,r_2)$, where $r_1\!:= -\frac{2}{3}( \ka - \frac{1}{2}c)^2$ and $r_2:= 2\ka c +  \frac{1}{2} c^2$.  \obsc{equivalencia} elucidates the fact that for such $r$, the set of smooth periodic TWS forms a continuous family 
$\{\vp_a\}_{a\in (0,\Am)}$ parameterized by the wave height as a consequence of \lemc{leng}. Moreover,  the wave length $\lambda(a)$ of the smooth periodic TWS $\vp_a$ is equal to the period of the periodic orbit of~\refc{e-sys_w} inside the energy level $\{\hat{H}=\ell^{-1}(a)\}.$ Hence, by applying \teoc{res_periode}, the result will follow once we write the conditions $\be \in\left(0,-\frac{1}{10}+\frac{1}{15}\sqrt{6}\,\right]$, $\be \in\left(-\frac{1}{10}+\frac{1}{15}\sqrt{6},\frac{1}{6}\right)$ and $\be \geqslant\frac{1}{6}$ in terms of $\ka$, $c$ and $r$. Taking the relation $\be=\frac{1}{6}\!\left(\frac{2|\beta |}{\sqrt{4\beta ^2+6\alpha }}-1\right)$ into account and setting 
 \[ 
 r_{b_1}\!:=\ka c-\frac{1}{2} \ka^2\,\text{ and }r_{b_2}\!:=\frac{\sqrt{6}-3}{6}\bigl((\sqrt{6}+1)\ka^2-2(\sqrt{6}-5)\ka c-2c^2\bigr),
 \] some computations show that these conditions are given, respectively, by $r\in [r_{b_2},r_2),$ $r\in (r_{b_1},r_{b_2})$ and $r\in (r_1,r_{b_1}).$ This proves the result.  
\end{prooftext}

%%%%%%%%%%%%%%%

\section{Study of the period function}\label{s-mon}

This Section is devoted to the proof of  \teoc{res_periode}, which strongly relies on the tools developed in \cite{Vil2014}. In order to explain how they can be applied to our problem, some definitions need to be introduced. 
 In the aforementioned paper the authors consider analytic planar differential systems
\begin{equation}\label{eq1}
 \sist{p(x,y),}{q(x,y),}
\end{equation}
satisfying the following \emph{hypothesis}:
\begin{center}
\textbf{(H)}\hspace{0.5truecm}
\begin{minipage}{0.75\linewidth}
The differential system \refc{eq1} has a center at the origin and an analytic first integral of the form $H(x,y)=A(x)+B(x)y+C(x)y^2$ with $A(0)=0$. Moreover its integrating factor, say $\infa$, depends only on~$x.$
\end{minipage}\phantom{\textbf{(H)}\hspace{0.5truecm}}
\end{center}
Let $\ii$ be the projection onto the $x$-axis of the period annulus $\P$ around the center at the origin of the differential system \refc{eq1}. Note that $x_{\ell}<0<x_r.$ Then, by Lemma 3.1 in~\cite{Vil2014}, the hypothesis \textbf{(H)} implies that $M\!:=\frac{4AC-B^2}{4|C|}$ is a well defined analytic function on $\ii$ with $M(0)=0$ and $xM'(x)>0$ for all $x\in\ii\setminus\{0\}.$ Accordingly, there exists a unique analytic function $\sigma$ on $\ii$ with $\sigma(x)=-x+\op(x)$ such that $M\circ\sigma=M.$ Note that $\sigma$ is an \emph{involution} with $\sigma(0)=0.$ (Recall that a mapping $\sigma$ is said to be an involution if $\sigma\circ\sigma=\mbox{Id}.)$ Given an analytic function $f$ on $\ii\setminus\{0\}$ we define its $\sigma$-\emph{balance} to be
\[
 \mathscr B_{\sigma}\bigl(f\bigr)(x)\!:=\frac{f(x)-f\bigl(\sigma(x)\bigr)}{2}.
\]
Taking these definitions into account, the statement $(b)$ of \teoc{thm}  in \cite{Vil2014} asserts the following:

\begin{prop}\label{l-Vil14}
Suppose that the analytic differential system \refc{eq1} satisfies the hypothesis~\textbf{(H)}. Setting $\mu_0=-1$, define recursively
\[
 \mu_{i}\!:=\left(\frac{1}{2}+\frac{1}{2i-3}\right)\mu_{i-1}+\frac{\sqrt{|C|}M}{(2i-3)\infa}
  \left(\frac{\infa\mu_{i-1}}{\sqrt{|C|}M'}\right)'
  \mbox{ and }\;
 \ell_i\!:=\frac{\infa\mu_i}{\sqrt{|C|}M'}\,
 \mbox{ for $i\geqslant 1.$}
\]
If the number of zeros of $\mathscr B_{\sigma}(\ell_i)$ on $(0,x_r)$, counted with 
multiplicities, is $n\geqslant 0$ and it holds that $i>n$, then the number of critical periods of the center at the origin, 
 counted with multiplicities, is at most $n.$
 \end{prop}
In particular, we point out that the period function is monotonous if $n=0$. A key ingredient for determining the number of zeros of $ \mathscr B_{\sigma}\bigl(l_i\bigr)$  is the following result, see \cite[Theorem B]{Vil2014}. In its statement, and in what follows, $\emph{\mbox{Res}}$ stands for the \emph{multipolynomial resultant} (see for instance~\cite{Cox2007,Fulton1984}).

\begin{prop}\label{cota}
Let $\sigma$ be an analytic involution on $\ii$ with $\sigma(0)=0$ and let $\ell$ be an analytic function on $\ii\setminus\{0\}.$ Assume that $\ell$ and $\sigma$ are algebraic, i.e., that there exist $L,S\in\C[x,y]$ such that $L\bigl(x,\ell(x)\bigr)\equiv 0$ and $S\bigl(x,\sigma(x)\bigr)\equiv 0$. Let us define $T(x,y)\!:=\mbox{Res}_z\bigl(L(x,z),L(y,z)\bigr)$ and $\mathscr R(x)\!:={Res}_y\bigl(S(x,y),T(x,y)\bigr).$ Finally let $s(x)$ and $t(x)$ be, respectively, the leading coefficients of $S(x,y)$ and $T(x,y)$ with respect to $y$. Then the following hold: 
\begin{enumerate}[$(a)$]
\item If $\mathscr B_{\sigma}\bigl(\ell\bigr)(x_0)=0$ for some $x_0\in\ii\setminus\{0\},$ then $\mathscr R(x_0)=0.$ 
\item If $s(x)$ and $t(x)$ do not vanish simultaneously at $x_0,$ then the multiplicity of 
         $\mathscr B_{\sigma}\bigl(\ell\bigr)$ at $x_0$ is not greater than the multiplicity of $\mathscr R$ at $x_0.$ 
\end{enumerate}
\end{prop}

In order to apply these results we move the center of differential system \eqref{e-sys_w} to the origin. In passing we notice that the problem is essentially one-parametric. 
Since its proof is a straightforward computation, we do not include it here for the sake of brevity. 

\begin{lem}\label{l-lem1}
Consider system \refc{e-sys_w} with $\alpha $ and $\beta $ verifying $-2\beta ^2 < 3\alpha  < 0$ and let us say that the center is at a point $(w_c,0)$. Then the coordinate transformation given by $\left\{x=\frac{w-w_c}{2\beta \sqrt{\De}},y=\frac{v}{2\beta \sqrt{\De}}\right\}$, where $\De\!:= 4+\frac{6\alpha }{\beta ^2}$, brings system \refc{e-sys_w} to
\begin{equation}\label{e-sys_x}
       \left\{
      \begin{array}{l}
     x' =y,\\[2pt]
     y' = -\dfrac{x-3x^2 +y^2}{2(x+\be)},
      \end{array}\right.
\end{equation}
where $\be\!:= \frac{1}{6}\!\left(\frac{2}{\sqrt{\De}}-1\right)$ is positive.
\end{lem}
The planar differential system~\eqref{e-sys_x} is analytic away from the singular line $x=-\be .$ One can easily verify that it satisfies the hypothesis  \textbf{(H)} with $A(x)=\frac{1}{2}x^2-x^3$, $B(x)=0$, $C(x)=x+\be $ and $\infa(x)=2(x+\be ).$ The function $A$ has a minimum at $x=0$ and a maximum at $x=\frac{1}{3},$ which yield a center at $(0,0)$ and a saddle at $(\frac{1}{3},0),$ respectively. 
\begin{figure}[t]
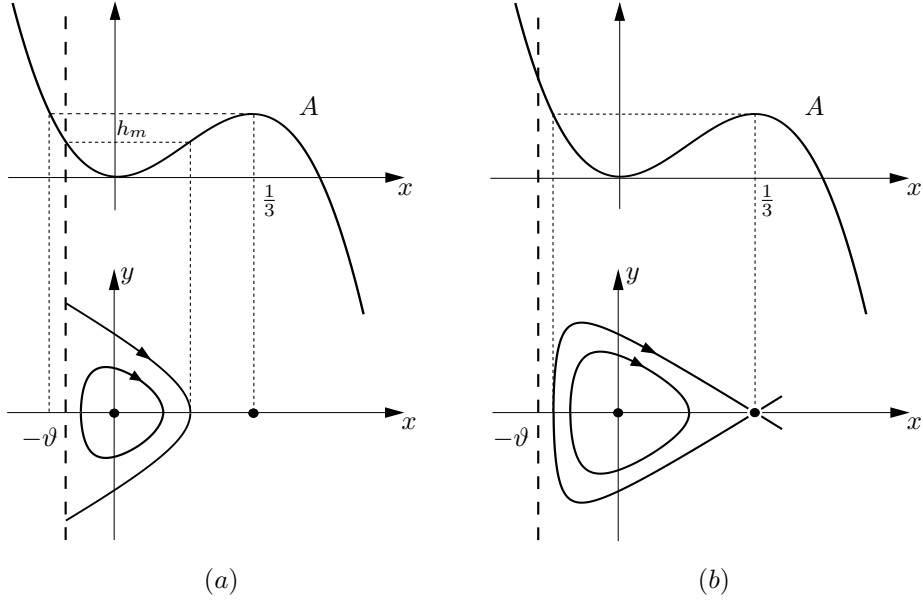

  \centering    
 \begin{lpic}[l(-10mm),r(2mm),t(0mm),b(0mm)]{Figure2(0.5)}
     \lbl[l]{131,108;$x$}
     \lbl[l]{105,130;$A$}
     \lbl[l]{132,46;$x$}
     \lbl[l]{58,85;$y$}

     \lbl[l]{32,42;$-\be $}
     \lbl[l]{95,105;$\frac{1}{3}$}
     \lbl[l]{57,124;\footnotesize{$h_m$}}
     
    \lbl[l]{265,108;$x$}
     \lbl[l]{238,130;$A$}
     \lbl[l]{265,46;$x$}
     \lbl[l]{191,85;$y$}
    
     \lbl[l]{156,42;$-\be $}
     \lbl[l]{226,105;$\frac{1}{3}$}
          
      \lbl[l]{80,4;$(a)$}
      \lbl[l]{210,4;$(b)$}
 \end{lpic}     
\caption{The period annulus $\P$ of the center at the origin of system \eqref{e-sys_x} for the two different cases that may occur: $(a)$ when $\be <1/ 6$; $(b)$ when $\be > 1/6$.}
\label{Fig-period annulus}
\end{figure}
When $\be>\frac{1}{6}$, in which case $A(-\be )>A(\frac{1}{3})$,  the singular line is ``far away'' from~$\P$, and the period annulus is bounded by the homoclinic connection based in the saddle point. When $\be<\frac{1}{6}$ the situation is quite different because the outer boundary of $\mathscr P$ consists of a trajectory with $\alpha$ and $\omega$ limit in the straight line $\{x=-\be \}$ and the segment between these two limit points, see~\figc{Fig-period annulus}. For this reason, we will study the period function of the center of system \eqref{e-sys_x} separately for $\be <1/6$ and $\be >1/6.$

Observe that if $B=0$, then the hypothesis \textbf{(H)}  implies that the involution~$\sigma$ is defined by $A=A\circ\sigma$. This is the case in the differential system under consideration, and one can easily verify that
 \begin{equation}\label{e-sig}
  A(x)-A(z)=2(z-x)S(x,z),\,\mbox{ where $S(x,z)\!:=2x^2 +2xz+2z^2-x-z$}.
 \end{equation}
Thus, we get $\sigma(x)=\frac{1}{4}(1-2x-\sqrt{(6x+1)(1-2x)}$. As a matter of fact, thanks to \propc{cota}, the explicit expression of the involution is not required and we shall only use that  $S\bigl(x,\sigma(x)\bigr)=0$.

The following auxiliary result will be needed at various points throughout this Section. The proof is a straightforward computation of the first three coefficients in the Taylor expansion of the period function using standard techniques (see for example \cite{Gasull1997}).
 
\begin{lem}\label{L-T0}
The first, second and third period constants of the center at the origin of system~\eqref{e-sys_x} are given, up to a positive factor, by 
 \[
  \Delta_1=60\be^2+12\be-1,\;\Delta_2=-\Delta_1\mbox{ and }\,
  \Delta_3=18240\,{\be }^{4}+3312\,{\be }^{3}-276\,{\be }^{2}+40\,\be -5,
 \]
respectively.
\end{lem}

\begin{prop}\label{prop: beta > 1/6}
If $\be \geqslant \frac{1}{6},$ then the period function of the center of system~\eqref{e-sys_x} is monotonous increasing.
\end{prop}

\begin{prova}
If $\be\geqslant\frac{1}6$ then, see \figc{Fig-period annulus}, the projection of the period annulus on the $x$-axis is $\bigl(-\frac{1}{6},\frac{1}{3}\bigr)$. Following \propc{l-Vil14}, we shall study the number of zeros of $\Bs (\ell_1)$ and to this end we will apply \propc{cota}. With this aim in view note that
\begin{equation*}
     \ell_1(x)=\frac{1}{2}\frac{(6\be+1)x-4\be-1}{\sqrt{x+\be}(3x-1)^3}.
\end{equation*}
Accordingly, $L\bigl(x,\ell_1(x)\bigr)\equiv 0$ with $L(x,y)\!:=4(x+\be )(3x-1)^6y^2-\bigl((6\be+1)x-4\be-1\bigr)^2.$ Recall also that $S\bigl(x,\sigma(x)\bigr)\equiv 0$, where $S\in\R[x,y]$ is defined in~\refc{e-sig}. A computation shows that $\mbox{Res}_z\bigl(L(x,z),L(y,z)\bigr)=16(x-y)^2\hat T(x,y)^2$, with $\hat T$ a bivariate polynomial of degree $8$ in~$x$ and~$y$ which also depends polynomially on~$\be$. Finally $\mathscr R(x)\!:=\mbox{Res}_y\bigl(S(x,y),T(x,y)\bigr)=(3x-1)^8 R(x)$, where~$R$ is a univariate polynomial of degree $8$ in $x$ depending polynomially on~$\be$. 

Let us define $\mathcal Z(\be )$ to be the number of roots of $R$ on $(0,\frac{1}{3})$ counted with multiplicities. We claim that $\mathcal Z(\be )=0$ for all $\be \geqslant\frac{1}{6}.$ For $\be =\frac{1}{6}$ this can be easily verified by applying Sturm's Theorem. To prove it for $\be >\frac{1}{6}$ we first note that
\begin{align}
  &R(0)=(4\be+1)(2\be+1)(48\be^2 + 24\be -1)(60\be^2 + 12 \be - 1)\label{eq5} \\
\intertext{and}
  &R(1/3)=(6\be-1)(2160\be^3+2484\be^2+720\be +17)(2/3+2\be)^2\notag, 
\end{align}
which do not vanish for $\be>\frac{1}{6}$. The discriminant of $R$ with respect to~$x$, $\mbox{Disc}_x(R)$, is a polynomial $\mathscr D(\be)$ of degree 82. After  factorizing it, one can easily prove that $\mathscr D$ vanishes on $(\frac{1}{6},+\infty)$ exactly once, at $\be =\bar\be $ with $\bar\be \approx 0.954.$ Altogether this implies that $\mathcal Z(\be )$ is constant on $(\frac{1}{6},\bar\be )$ and $(\bar\be ,+\infty).$ Choosing one value of~$\be$ in each interval and applying Sturm's Theorem we find that $\mathcal Z(\be )=0$ for $\be \in (\frac{1}{6},+\infty)\setminus\{\bar\be \}.$ To prove that this is true for $\be =\bb$ as well we show that $x\longmapsto \ell_1(x)$ is monotonous on $(-\frac{1}{6},\frac{1}{3})$ for all $\be\in(\frac{9}{10},1)$. Indeed, one can verify that 
\[
    \ell_1'(x)=\frac{N(x)}{4(x+\be)^{3/2}(3x-1)^4},
\]
with
\[
   N(x) = (90\be+15)x^2+(72\be^2-66\be-20)x -60\be^2-12\be+1.
\]
We have that $ N(x) \neq 0$ for  $x\in(-\frac{1}{6},\frac{1}{3})$ and $\be\in(\frac{9}{10},1)$ because it is true for $\be =\frac{95}{100}$ and, on the other hand, the number of roots counted with multiplicity does not change due to the fact that
\[
    N(-1/6) N(1/3)\,\mbox{Disc}_x(N)\neq 0 \mbox{ for all  }\be\in(9/10,1).
\]
Therefore $\ell_1'(x)\neq 0$ for all $x\in(-\frac{1}{6},\frac{1}{3})$ and $\be\in(\frac{9}{10},1)$,  and we can assert that $\mathscr R$ does not vanish on $(0,\frac{1}{3})$ for any $\be\geqslant\frac{1}{6}.$ In view of $(a)$ in \propc{cota} this implies that $\Bs\bigl(\ell_1\bigr)\neq 0$ on $(0,\frac{1}{3}).$ This proves the validity of the claim and hence, by applying \propc{l-Vil14} with $n=0,$ it follows that the period function is monotonous for $\be\geqslant\frac{1}{6}.$ Finally, the result follows by noting that, thanks to \lemc{L-T0}, the first period constant $\Delta_1$ is positive for $\be\geqslant \frac{1}{6}$. 
\end{prova}

In order to study the period function  of the center of system \eqref{e-sys_x} for $\be  < \frac{1}{6}$,  we first recall the well-known \emph{Gelfand-Leray derivative}, see for instance \cite{Ilyashenko2006}.

\begin{lem}\label{L-GL}
Let $\omega$ and $\eta$ be two rational $1$-forms such that $d \omega  = dH \wedge \eta$ and let  $\gamma_h \in H^1(L_h,\Z)$ be a continuous family of cycles on non-critical level curves $L_h=\{H=h\}$ not passing through poles of neither $\omega$ nor $\eta$. Then 
\begin{equation}
\label{eq-GF}
    \frac{d}{dh} \oint \omega = \oint \eta. 
\end{equation}
\end{lem}

We shall also use the following result, see \cite[Lemma 4.1]{Grau2011}.

\begin{lem}\label{L-GraVil14}
Let $\g_h$ be an oval inside the level curve $\{A(x)+C(x)y^2=h\}$ and consider a function $F$ such that $F/A'$ is analytic at $x=0$. Then, for any $k\in \N$, 
 \begin{equation*}
     \int_{\g_h}\!F(x) y^{k-2} dx = \int_{\g_h}\!G(x) y^k dx, 
 \end{equation*}
 where $G= \frac{2}{k}\left(\frac{CF}{A'}\right)'- \bigl(\frac{C'F}{A'}\bigr).$
\end{lem}
This allows us to prove the following result about the derivative of the period function associated to the center of the analytic differential system \refc{eq1} satisfying hypothesis \textbf{(H)}. 
\begin{lem}\label{T'}
Suppose that the analytic differential system \refc{eq1} satisfies the hypothesis \textbf{(H)} with $B=0.$ Let $T(h)$ be the period of the periodic  orbit $\gamma_h$ inside the energy level $\{H=h\}$. Then
 \[
  T'(h)=\frac{1}{h}\int_{\gamma_h}\!R(x)\,\frac{dx}{y}, 
 \]
where $R\!=\frac{1}{2C}\!\left(\frac{\infa A}{A'}\right)'-\frac{\infa(AC)'}{4A'C^2}$.
\end{lem}

\begin{prova}
Note first that if \refc{eq1} satisfies \textbf{(H)} with $B=0,$ then $\frac{dx}{dt}=\frac{H_y(x,y)}{\infa(x)}=\frac{2C(x)y}{\infa(x)},$ so that 
 \begin{equation}\label{eq2}
  T(h)=\int_{\gamma_h}\!\left(\frac{\infa}{2C}\right)\!(x)\frac{dx}{y}.
 \end{equation}
Accordingly, since $A(x)+C(x)y^2=h$ on $\gamma_h$ we get
 \[
  2hT(h)=\int_{\gamma_h}\!\left(\frac{\infa A}{C}\right)\!(x)\,\frac{dx}{y}+\int_{\gamma_h}\!\infa(x)ydx=
    \int_{\gamma_h}\!\bigl(G+\infa\bigr)(x)\,ydx,
 \]
with $G\!:=2\left(\frac{\infa A}{A'}\right)'-\frac{\infa AC'}{A'C}$, where the second equality follows by applying \lemc{L-GraVil14} with $F=\frac{\infa A}{C}.$ Next we apply \lemc{L-GL} taking $H(x,y)=A(x)+C(x)y^2$, $\omega=\bigl(G+\infa\bigr)(x)\,ydx$ and $\eta=\left(\frac{G+\infa}{2C}\right)\!(x)\frac{dx}{y}$ in order to get that
 \[
  2\bigl(hT(h)\bigr)'=2hT'(h)+2T(h)=\int_{\gamma_h}\!\left(\frac{G+\infa}{2C}\right)\!(x)\frac{dx}{y}.
 \]
This equality, on account of \refc{eq2}, implies that $2hT'(h)=\int_{\gamma_h}\!\left(\frac{G-\infa}{2C}\right)\!(x)\frac{dx}{y}$. This proves the result because a straightforward computation shows that $R=\frac{G-\infa}{4C}$. 
\end{prova}

We are now in position to prove  the following:

\begin{lem}
\label{L-T'lim}
If $\be  \in (0,\frac{1}{6})$, then the period function $T(h)$ of the center at the origin of~\eqref{e-sys_x} verifies that $\lim_{h\to h_m} T'(h)=-\infty$, where $h_m =A(-\be )$ is the energy level of the outer boundary of $\P$, see~\figc{Fig-sketch}.
\end{lem}

\begin{prova}
By applying \lemc{T'} taking $A(x)=\frac{1}{2}x^2-x^3$, $C(x)=x+\be $ and $\infa(x)=2(x+\be )$  it follows that $T'(h)=\frac{1}{h}\int_{\gamma_h}\!R(x)\frac{dx}{y}$ with
 \[
  R(x)\!:=\frac{x\bigl(4\be+1-(6\be+1)x\bigr)}{ 4(x+\be)(3x-1)^2}.
 \]
The relative position of the straight line $x=-\be$ with respect to the graph of $A$ is as  displayed in \figc{Fig-sketch}  because $A(-\be)-A(\frac{1}{3})=\frac{1}{54}(6\be-1)(3\be+1)^2$ and, by assumption, $\be\in (0,\frac{1}{6}).$ Accordingly if $h\in (0,\frac{1}{54})$ then $h-A(x)=(x-x_h^-)(x-x_h^+)(x-\hat x_h)$, where $x_h^-<0<x_h^+<\frac{1}{3}<\hat x_h$. In particular, for $h\in (0,h_m)$, the projection of the periodic orbit $\gamma_h$ on the $x$-axis is the interval $[x_h^-,x_h^+]$. 
\begin{figure}[t]
  \centering    
 \begin{lpic}[l(0mm),r(0mm),t(-10mm),b(-30mm)]{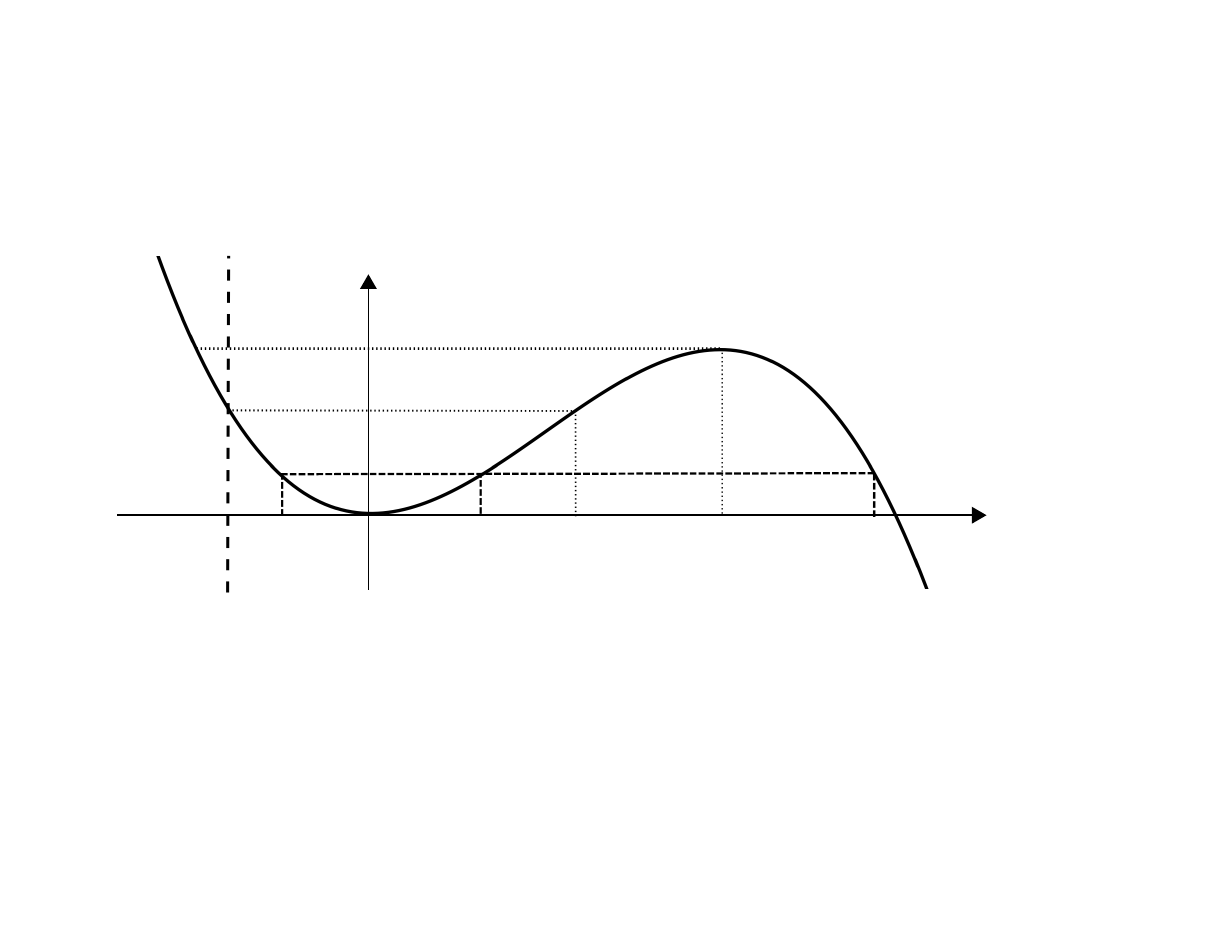}
     \lbl[l]{99,39;$x$}
     \lbl[l]{83,58;$A$}
     \lbl[l]{17,39;$-\be $}
     \lbl[l]{27,39;$x_h^-$}
     \lbl[l]{47,39;$x_h^+$}
     \lbl[l]{57,39;$x_r$}
     \lbl[l]{72,39;$\frac{1}{3}$}
     \lbl[l]{87.5,39;$\hat{x}_h$}
     \lbl[l]{38.5,48.5;$h$}
     \lbl[l]{38,55;$h_m$}
     \lbl[l]{38,62;$\frac{1}{54}$}
 \end{lpic}     
\caption{Root distribution of $h-A(x)=(x-x_h^-)(x-x_h^+)(x-\hat x_h)$ in the proof of \lemc{L-T'lim}}\label{Fig-sketch}
\end{figure}
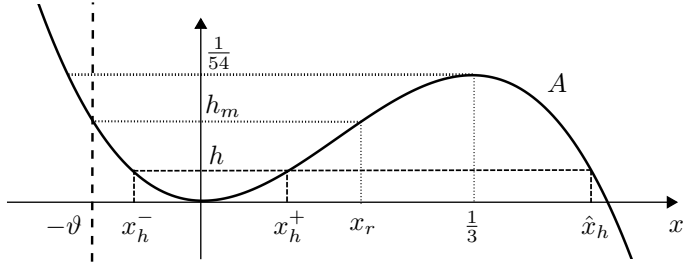
Hence $T'(h)=\frac{2}{h}\bigl(I_1(h)+I_2(h)\bigr),$ where
 \[ 
  I_1(h)=\int_{x_h^-}^0f(x,h)\,dx\,\text{ and }
  I_2(h)=\int_0^{x_h^+}f(x,h)\,dx
 \]
 with
 \[
  f(x,h)=\frac{R(x)\sqrt{C(x)}}{\sqrt{h-A(x)}}=\frac{x\bigl(4\be+1-(6\be+1)x\bigr)}{ 4(3x-1)^2\sqrt{x+\be}
          \sqrt{(x-x_h^-)(x-x_h^+)(x-\hat x_h)}}.
 \]
Let us write $f(x,h)=\frac{g_1(x,h)}{\sqrt{(x+\be )(x-x_h^-)}}$, where 
 \[ 
  g_1(x,h)\!:=\frac{x\bigl(4\be+1-(6\be+1)x\bigr)}{ 4(3x-1)^2
          \sqrt{(x-x_h^+)(x-\hat x_h)}}.
 \]  
Note that $g_1$ is a continuous function on $(-\infty,0]\!\times\! (0,\frac{1}{54}).$ Consequently there exists $M_1\in\R$ such that $M_1\!:=\sup\bigl\{g_1(x,h);(x,h)\in [-\frac{1}6,0]\!\times\! [\frac{1}{2}h_m,h_m]\bigr\}.$ In addition, observe that $M_1$ is strictly negative because one can verify that $4\be+1-(6\be+1)x>0$ for all $x<0$ and $\be>0.$ Thus for $h\in (\frac{1}{2}h_m,h_m)$ we have that
  \begin{align*}
     I_1 (h)&= \int_{x_h^-}^0\frac{g_1(x,h)dx}{\sqrt{(x+\be)(x-x_h^-)}}
            \leqslant M_1 \int_{x_h^-}^0\frac{dx}{\sqrt{(x+\be)(x-x_h^-)}} \\
          &= M_1 \log\left(\frac{\be-x_h^- +\sqrt{-\be x_h^-}}{\be+x_h^-}\right)
             \longrightarrow -\infty \,\text{ as } h\longrightarrow h_m.
 \end{align*}
In the inequality above we take $-\frac{1}{6}<-\be<x_h^-$ into account, whereas the limit follows by using $ M_1<0$ and the fact that $x_h^-$ tends to $-\be$ as $h\longrightarrow h_m.$ Accordingly, 
 \begin{equation}\label{eq4}
 \lim_{h\to h_m} I_1(h)=-\infty.
\end{equation}
In order to study $I_2$ let us write $f(x,h)=\frac{g_2(x,h)}{\sqrt{x_h^+-x}}$, where 
 \[ 
  g_2(x,h)\!:=\frac{x\bigl(4\be+1-(6\be+1)x\bigr)}{ 4(3x-1)^2
          \sqrt{(x+\be)(x-x_h^-)(\hat x_h-x)}}.
 \]  
Since $g_2$ is continuous on $[0,\frac{1}{3})\!\times\!(0,\frac{1}{54}),$  $M_2\!:=\sup\bigl\{g_2(x,h);(x,h)\in [0,x_r]\!\times\! [\frac{1}{2}h_m,h_m]\bigr\}$ is a well defined real number. Consequently if $h\in (\frac{1}{2}h_m,h_m),$ then
  \[
     I_2 (h)= \int^{x_h^+}_0\frac{g_2(x,h)dx}{\sqrt{x_h^+-x}}
            \leqslant M_2\int^{x_h^+}_0\frac{dx}{\sqrt{x_h^+-x}} 
            =2M_2\sqrt{x_h^+}<\frac{2M_2}{\sqrt{3}}.
 \]
Due to $T'(h)=\frac{2}{h}\bigl(I_1(h)+I_2(h)\bigr)$, the above inequality together with \refc{eq4} imply the result.    
\end{prova}

\begin{prooftext}{Proof of \teoc{res_periode}.}
Thanks to \propc{prop: beta > 1/6} it suffices to consider $\be \in (0,\frac{1}{6})$. For these parameter values, see \figc{Fig-sketch}, the projection of the period annulus on the $x$-axis is $(-\be ,x_r)$, where $A(x_r)=A(-\be)$. We proceed in exactly the same way as we did with \propc{prop: beta > 1/6}, i.e., by applying \propc{l-Vil14} together with \propc{cota}, but in this case we must use $\ell_3,$ since neither $\ell_1$ nor $\ell_2$ provide decisive information. Since $\mathscr B_{\sigma}(f)\circ\sigma=-\mathscr B_{\sigma}(f)$ and $\sigma$ maps $(0,x_r)$ to $(x_{\ell},0),$ for convenience we shall study the latter interval, which in this case is $(-\be ,0).$ One can verify that
 \[
  \ell_3(x)=\frac{p(x)}{(3x-1)^7(x+\be)^{5/2}},
 \]
where $p$ is a polynomial of degree 7 in $x$ (depending also polynomially on $\be $), which we do not write for the sake of  brevity. Therefore $L\bigl(x,\ell_3(x)\bigr)\equiv 0$ with $L(x,y)\!:=(x+\be )^5(3x-1)^{14}y^2-p(x)^2.$ Recall also that $S\bigl(x,\sigma(x)\bigr)\equiv 0$, where $S$ is the polynomial given in~\refc{e-sig}. A computation shows that $\mbox{Res}_z\bigl(L(x,z),L(y,z)\bigr)=2^{-20}(x-y)^2\hat T(x,y)^2$, with $\hat T\in\R[x,y]$ of degree $32$, depending also polynomially on~$\be$. Finally $\mathscr R(x)\!:=\mbox{Res}_y\bigl(S(x,y),T(x,y)\bigr)=(3x-1)^{20} R(x)$, where~$R\in\R[x,\be]$ with $\deg(R;x)=44$. For each $\be\in (0,\frac{1}{6})$ let us define $\mathcal Z(\be)$ to be the number of zeros, counted with multiplicities, of $R$ on $(-\be,0).$ To study this number we consider the value of~$R$ at the endpoints of $(-\be,0)$,
 \begin{align*}
  &R(0)=2^{11}5^3\,{\be }^{12} \left( 1+4\,\be  \right)  \left( 60\,{\be }^{2}
+12\,\be -1 \right)  \left( 48\,{\be }^{2}+24\,\be -1 \right) 
 \left( 2\,\be +1 \right) ^{5},\\
 \intertext{and}
 &R(-\be)=16\,{\be }^{12} \left( 1+3\,\be  \right) ^{10} \left( 2\,\be +1
 \right) ^{12} \left( 6\,\be -1 \right) ^{14},
 \end{align*}
together with the discriminant of $R$ with respect to $x,$ $\mbox{Disc}_{x}(R),$ which is a polynomial~$\mathscr D(\be)$ of degree 1586 that we do not write here for brevity. 
 One can easily check that $R(-\be )$ does not vanish and that $R(0)$ has exactly two roots on $(0,\frac{1}{6})$, namely 
 \[
  \be_1\!:=-\frac{1}{4}+\frac{1}{6}\sqrt{3}\approx 0.03867\,\mbox{ and }
  \be_2\!:=-\frac{1}{10}+\frac{1}{15}\sqrt{6}\approx 0.06330
 \]
By applying Sturm's Theorem to each factor, we conclude that on  $(0,\frac{1}{6})$ the discriminant $\mathscr D(\be)$ vanishes only at $\be=\be_2$. 
Hence $\mathcal Z(\be)$ is constant on $I_1\!:=(0,\be_1),$ $I_2\!:=(\be_1,\be_2)$ and $I_3\!:=(\be_2,\frac{1}{6}).$ Taking one parameter value on each interval and applying Sturm's Theorem once again we can assert that $\mathcal Z(\be)=0$ for all $\be\in I_1,$ $\mathcal Z(\be)=1$ for all $\be\in I_2$ and $\mathcal Z(\be)=2$ for all $\be\in I_3.$ Therefore, by \propc{cota}, it follows that the number of zeros, counted with multiplicities, of $\mathscr B_{\sigma}\bigl(\ell_3\bigr)(x)$ on $(-\be,0)$ is at most 0, 1 and 2, for $\be\in I_1$, $\be\in I_2$ and $\be\in I_3$, respectively. Hence, thanks to \propc{l-Vil14}, we can assert that the period function is monotonous for $\be\in I_1$, whereas it has at most 1 (respectively, 2) critical periods for $\be\in I_2$ (respectively, $\be\in I_3$), counted with multiplicities.  
 
Recall at this point that, in view of \lemc{L-T0}, the first period constant of the center is given by $\Delta_1=60\be^2+12\be-1$. On the other hand, by \lemc{T'}, we know that $\lim T'(h)=-\infty$ as $h$ tends to $h_m$ for all $\be \in (0,\frac{1}{6}).$ Since $\Delta_1=0$ for $\be=-\frac{1}{10}\pm\frac{1}{15}\sqrt{6},$ we conclude that $T(h)$ is monotonous decreasing near the endpoints of $(0,h_m)$ for all $\be \in (0,\be _2).$ For the same reason, if $\be\in I_3$ then $T(h)$ is increasing near $h=0$ and decreasing near $h=h_m.$ On account of the upper bounds on the number of critical  periods that we have previously obtained, we conclude that the period function is monotonous decreasing for $\be\in I_1\cup I_2$  and it has a unique critical period, which is a maximum, for $\be\in I_3$.

The fact that the period function is monotonous decreasing for $\be=\be_1$ can be proved by showing that $\mathscr B_{\sigma}\bigl(\ell_1\bigr)$ does not vanish on $(-\be_1,0)$ and using that $\Delta_1<0$ at $\be=\be_1.$ Since this is easy we do not include it here for the sake of brevity. The proof for $\be=\be_2$ is slightly different but straightforward as well. We show first that $\mathscr B_{\sigma}\bigl(\ell_3\bigr)$ has at most one zero on $(-\be_2,0)$ counted with multiplicities.  By \propc{l-Vil14} this implies that the period function has at most one critical period. To prove that it has none we take the behaviour of the period function at the endpoints of $(0,h_m)$ into account. Since $\Delta_1=\Delta_2=0$ and $\Delta_3<0$ at $\be=\be_2$ by \lemc{L-T0}, we have that it is decreasing near $h=0.$ We know that it is also decreasing near $h=h_m$ thanks to \lemc{T'}. Thus it can not have any critical period. This completes the proof.
\end{prooftext} 
 
 \begin{obs}\label{remark}
By means of standard techniques one can obtain the limit value of the integral defining the period function at the endpoints of its interval of definition. 
Combining this information with the results in \teoc{res_periode} and \secc{s-mon} we get the graphs of the period function $T(h)$ displayed in Figure~\ref{Fig-periodfunction}.
\begin{figure}[t]
  \centering    
 \begin{lpic}[l(-2mm),r(0mm),t(5mm),b(10mm)]{Figure4(0.40)}
      \lbl[l]{45,-15;(a)}
      \lbl[l]{8,95;$T$}
      \lbl[l]{110,5;$h$}      
      \lbl[l]{15,-2;$h_0$}
      \lbl[l]{84,-2;$h_1$}
      \lbl[c]{0,68;$T_0$}
      \lbl[c]{0,23;$T_1$}

      \lbl[l]{175,-15;(b)}
      \lbl[l]{130,95;$T$}
      \lbl[l]{232,5;$h$}  
      \lbl[l]{136,-2;$h_0$}
      \lbl[l]{207,-2;$h_1$}
      \lbl[c]{122,58;$T_0$}
      \lbl[c]{122,18;$T_1$}
            
      \lbl[l]{305,-15;(c)}        
      \lbl[l]{252,95;$T$}
      \lbl[l]{354,5;$h$}
      \lbl[l]{262,-2;$h_0$}
      \lbl[l]{334,-2;$h_1$}
      \lbl[c]{244,15;$T_0$}
 \end{lpic}     
\caption{Sketch of the graph of the period function $T(h)$ corresponding to \teoc{res_periode}: (a) for $\be \in(0,-\frac{1}{10}+\frac{1}{15}\sqrt{6}\,]$; (b) for $\be \in(-\frac{1}{10}+\frac{1}{15}\sqrt{6},1/6)$; and (c) for  $\be \in[1/6,\infty)$.
}
\label{Fig-periodfunction}
\end{figure}
For the sake of brevity we omit the computations of the explicit values 
\begin{equation*}
    T_0 = 2\pi \sqrt{2\be}\,\text{ and }T_1 =  2 \ln \left( \frac{{(2\be+1)(1-6\be)}}{1+6\be-4\sqrt{\be(1+3\be)}}\right).
\end{equation*}
Note that $T_0$ and $T_1$ are strictly positive whenever they are defined. Taking into account the relation between period and wave length, cf.~\obsc{equivalencia}, this shows that there do not exist smooth periodic TWS of CH with arbitrarily small wave length. 
Finally, we point out that $\be =-\frac{1}{10}+\frac{1}{15}\sqrt{6}$ and $\be=1/6$ are, respectively, bifurcation values of the period function at the inner and outer boundary of the period annulus, see~\cite{Mardesi2006}.
\end{obs}

%\bibliographystyle{abbrv}
%\bibliography{/Users/anna/Documents/UNI/library.bib}

\end{document}